\documentclass[]{amsart}

\usepackage{amssymb}
\usepackage{url}
\usepackage{amscd}
\usepackage[dvips]{color}
\usepackage{graphicx}
\usepackage{epsfig}
\usepackage{subfigure}

\theoremstyle{plain}
\newtheorem{theorem}{Theorem}

\newtheorem{corollary}[theorem]{Corollary}

\newtheorem{proposition}[theorem]{Proposition}
\newtheorem{definition}{Definition}

\theoremstyle{remark}
\newtheorem{remark}{Remark}
\newtheorem{example}{Example}

\def\R{{\mathbb R}}
\def\C{{\mathbb C}}
\def\D{{\mathbb D}}
\def\H{{\mathbb H}}
\def\E{{\mathbb E}}

\def\v{{\bf v}}

\def\u{{\bf u}}
\def\1{{\bf 1}}
\def\i{{\bf i}}
\def\j{{\bf j}}
\def\k{{\bf k}}
\def\d{\partial}

\def\dis{\displaystyle}
\def\tr{\mathop{\rm tr}\nolimits}

\allowdisplaybreaks

\input{epsf.tex}

\begin{document}

\date{}

\title[Timelike Minimal Surfaces in Minkowski $3$-Space]{Weierstra{\ss} Representation for Timelike Minimal
  Surfaces in Minkowski $3$-Space}

\author{Sungwook Lee}

\address{Department of Mathematics\\
University of Southern Mississippi\\
Hattiesburg, MS 39406-5045, U.S.A.}

\email{sunglee@usm.edu}



\begin{abstract}
Using techniques of integrable systems, we study a Weierstra{\ss}
representation formula for timelike surfaces with prescribed mean
curvature in Minkowski $3$-space. It is shown that timelike minimal
surfaces are obtained by integrating a pair of Lorentz holomorphic
and Lorentz antiholomorphic null curves in Minkowski $3$-space. The
relationship between timelike minimal surfaces and bosonic
Nambu-Goto string worldsheets in spacetime is also discussed in the appendix.
\end{abstract}

\maketitle

\section{Preliminaries}
In this section, we review some basics on the geometry of timelike
surfaces in Minkowski $3$-space.

Let $\E^4_2$ be the semi-Euclidean $4$-space with rectangular
coordinates $x_0$,$x_1$,$x_2$,\\$x_3$ and the semi-Riemannian metric
$\langle\ ,\ \rangle$ of signature $(-,-,+,+)$ given by the
quadratic form
$$ds^2=-dx_0^2-dx_1^2+dx_2^2+dx_3^2.$$
The semi-Euclidean $4$-space $\E^4_2$ is identified with the linear space
$M_2\R$ of all $2\times 2$ real matrices via the correspondence
\begin{equation}
\label{eq:matrix}
\u=(x_0,x_1,x_2,x_3)\longleftrightarrow\begin{pmatrix}
x_0+x_3 & x_1+x_2\\
-x_1+x_2 & x_0-x_3
\end{pmatrix}.
\end{equation}
This identification is an isometry, since
$$\langle\u,\v\rangle=\frac{1}{2}\{\tr(\u\v)-\tr(\u)\tr(\v)\},\ \u,\v\in M_2\R.$$
In particular, $\langle\u,\u\rangle=-\det\u$. The standard basis
$e_0,e_1,e_2,e_3$ for $\E^4_2$ is identified with with the matrices
$$\1=\begin{pmatrix}
1 & 0\\
0 & 1
\end{pmatrix},\ \i=\begin{pmatrix}
0 & 1\\
-1 & 0
\end{pmatrix},\ \j^{'}=\begin{pmatrix}
0 & 1\\
1 & 0
\end{pmatrix},\ \k^{'}=\begin{pmatrix}
1 & 0\\
0 & -1
\end{pmatrix},$$
i.e.
\begin{equation}
(x_0,x_1,x_2,x_3)\longleftrightarrow x_0\1+x_1\i+x_2\j^{'}+x_3\k^{'}.
\end{equation}
Note that the $2\times 2$ matrices $x_0\1+x_1\i+x_2\j^{'}+x_3\k^{'}$
form the algebra $\H^{'}$ of \emph{para-quaternions} or
\emph{split-quaternions}. The basis $\{\1,\i,\j^{'},\k^{'}\}$
satisfies the following relation:
\begin{align*}
\i^2&=-\1,\ {\j^{'}}^2={\k^{'}}^2=\1,\\
\i\j^{'}&=-\j^{'}\i=\k^{'},\ \j^{'}\k^{'}=-\k^{'}\j^{'}=-\i,\
\k^{'}\i=-\i\k^{'}=\j^{'}.
\end{align*}
The imaginary split-quaternions
$${\rm Im}\H^{'}=\{x_1\i+x_2\j^{'}+x_3\k^{'}: x_i\in\R, i=1,2,3\}$$
is identified with Minkowski $3$-space $\E^3_1$.

Let $M$ be a connected orientable $2$-manifold and $\varphi:
M\longrightarrow\E^3_1$ an immersion.The immersion $\varphi$ is said
to be \emph{timelike} if the induced metric $I$ on $M$ is
Lorentzian. The induced Lorentzian metric $I$ determines a Lorentz
conformal structure ${\mathcal C}_I$ on $M$. Let $(x,y)$ be a
Lorentz isothermal coordinate system with respect to the conformal
structure ${\mathcal C}_I$. Then the first fundamental form $I$ is
written in terms of $(x,y)$ as
$$I=e^\omega(-dx^2+dy^2).$$
$(M,C_I)$ is said to be a \emph{Lorentz surface}.

Let $N$ be a unit normal vector field of $M$. Then $$ \langle
N,N\rangle=1,\
\langle\varphi_x,N\rangle=\langle\varphi_y,N\rangle=0.$$ The
\emph{conformality conditions} are given by
\begin{equation}
\label{eq:conf}
\left\{\begin{aligned} \langle\varphi_x,\varphi_y\rangle&=0,\\
-\langle\varphi_x,\varphi_x\rangle&=\langle\varphi_y,\varphi_y\rangle=e^\omega.
\end{aligned}\right.
\end{equation}

Let us define the quantities:
\begin{align*}
E&=\langle\varphi_x,\varphi_x\rangle,\
F=\langle\varphi_x,\varphi_y\rangle,\
G=\langle\varphi_y,\varphi_y\rangle,\\
l&=\langle\varphi_{xx},N\rangle,\ m=\langle\varphi_{xy},N\rangle,\
n=\langle\varphi_{yy},N\rangle.
\end{align*}
Then the mean curvature $H$ can be computed by the well-known
classical formula\footnote{This formula is still valid for spacelike
or timelike surfaces in Minkowski $3$-space.}
\begin{equation}
\label{eq:meancurv}
\begin{aligned}
H&=\frac{Gl+En-2Fm}{2(EG-F^2)}\\
 &=\frac{1}{2}e^{-\omega}\langle\Box\varphi,N\rangle,
\end{aligned}
\end{equation}
where $\Box$ is the \emph{d'Alembertian operator}
$$\Box=-\frac{\partial^2}{\partial x^2}+\frac{\partial^2}{\partial
y^2}.$$

Let $u:=x+y$ and $v=-x+y$. Then $(u,v)$ defines a null coordinate
system with respect to the conformal structure ${\mathcal C}_I$. The
first fundamental form $I$ is written in terms of $(u,v)$ as
$$I=e^\omega dudv.$$

The partial derivatives $\varphi_u$ and $\varphi_v$ are computed to
be
$$\varphi_u=\frac{1}{2}(\varphi_x+\varphi_y),\
\varphi_v=\frac{1}{2}(-\varphi_x+\varphi_y).$$ In terms of null
coordinates $(u,v)$,
$$\langle\varphi_u,N\rangle=\langle\varphi_v,N\rangle=0,$$
and the conformality condition \eqref{eq:conf} can be written as
\begin{equation}
\label{eq:conf2} \left\{\begin{aligned}
\langle\varphi_u,\varphi_u\rangle&=\langle\varphi_v,\varphi_v\rangle=0,\\
\langle\varphi_u,\varphi_v\rangle&=\frac{1}{2}e^\omega.
\end{aligned}\right.
\end{equation}
The mean curvature formula \eqref{eq:meancurv} is written as
\begin{equation}
\label{eq:meancurv2} H=2e^{-\omega}\langle\varphi_{uv},N\rangle.
\end{equation}

On a simply connected null coordinate region $\D$, we can find an orthonormal
frame field ${\mathcal F}: \D\longrightarrow{\rm O}^{++}(3,1)$ given by
\begin{align*}
{\mathcal F}&=(e^{-\omega/2}\varphi_x,e^{-\omega/2}\varphi_y,N)\\
            &=(e^{-\omega/2}(\varphi_u-\varphi_v),e^{-\omega/2}(\varphi_u+\varphi_v),N).
\end{align*}
Here, ${\rm O}^{++}(3,1)$ denotes the identity component of the Lorentz group
${\rm O}(3,1)$.

The special linear group ${\rm SL}(2,\R)$ acts isometrically on $\E^3_1$ via
the Ad-action:
$${\rm Ad}: {\rm SL}(2,\R)\times\E^3_1\longrightarrow\E^3_1;\ {\rm
  Ad}(g)X=gXg^{-1}.$$
The Ad-action induces a double covering ${\rm SL}(2,\R)\longrightarrow{\rm
  O}^{++}(3,1)$.

Using this double covering, we can find a lift $\Phi$ (called a
\emph{coordinate frame}) of ${\mathcal F}$ to ${\rm SL}(2,\C)$:
$${\rm Ad}(\Phi)(\i,\j^{'},\k^{'})={\mathcal F}.$$

${\mathfrak s}:=(\varphi_u,\varphi_v,N)$ defines a moving frame on $M$ and
satisfies the following Gau{\ss}-Weingarten equations:
\begin{equation}
\label{eq:g-w}
{\mathfrak s}_u={\mathfrak s}{\mathcal U},\ {\mathfrak s}_v={\mathfrak
  s}{\mathcal V},
\end{equation}
where
$${\mathcal U}=\begin{pmatrix}
\omega_u & 0 & -H\\
0 & 0 & -2e^{-\omega}Q\\
Q & \frac{1}{2}e^\omega H & 0
\end{pmatrix},\ {\mathcal V}=\begin{pmatrix}
0 & 0 & -2e^{-\omega}R\\
0 & \omega_v & -H\\
\frac{1}{2}e^\omega H & R & 0
\end{pmatrix}.$$
Here, $Q:=\langle\varphi_{uu},N\rangle$ and
$R:=\langle\varphi_{vv},N\rangle$. The quadratic differential
$1$-form
\begin{equation}
\label{eq:hopf} {\mathcal Q}=Qdu^2+Rdv^2
\end{equation}
is called \emph{Hopf differential}\footnote{This definition of Hopf
differential was suggested to the author by J. Inoguchi
\cite{inoguchi}.}. The second fundamental form $I\!I$ of $M$ is
related to ${\mathcal Q}$ by
\begin{equation}
\label{eq:sff} I\!I={\mathcal Q}+HI.
\end{equation}
This formula
implies that $p\in M$ is an \emph{umbilic point} if and only if
${\mathcal Q}(p)=0$, i.e., $p$ is a common zero of $Q$ and $R$.

If $K$ denotes the Gau{\ss}ian curvature, the Gau{\ss} equation
which describes a relationship between $K$, $H$, $Q$, and $R$ takes
the following form:
\begin{equation}
\label{eq:gauss2}
H^2-K=4e^{-2\omega}QR.
\end{equation}

The integrability condition for the Gau{\ss}-Weingarten equations is
\begin{equation}
\label{eq:g-c}
{\mathcal V}_u-{\mathcal U}_v+[{\mathcal U},{\mathcal V}]=0,
\end{equation}
which is equivalent to the Gau{\ss}-Codazzi equations
\begin{align}
\label{eq:gauss}
\omega_{uv}+\frac{1}{2}H^2e^\omega-2e^{-\omega}QR=0\\
\label{eq:codazzi}
\ H_u=2e^{-\omega}Q_v,\ H_v=2e^{-\omega}R_u.
\end{align}
It follows from the Codazzi equations \eqref{eq:codazzi} that mean curvature
$H$ is constant if and only if $Q_v=R_u=0$, i.e., $Q=Q(u)$ and $R=R(v)$. In
this case, $Q$ and $R$ are said to be \emph{Lorentz holomorphic} and
\emph{Lorentz anti-holomorphic}, respectively. For more details about Lorentz
holomophicity, see \cite{D-I-T} and \cite{I-T}.

Let $\Phi: \D\longrightarrow{\rm SL}(2,\R)$ be a coordinate frame of
$\varphi$. Then the partial derivatives $\varphi_u$ and $\varphi_v$ are given
in terms of the coordinate frame $\Phi$ as
\begin{align*}
\varphi_u&=\frac{1}{2}(\varphi_x+\varphi_y)\\
         &=\frac{1}{2}e^{\omega/2}{\rm Ad}(\i+\j^{'})\\
         &=e^{\omega/2}\Phi\begin{pmatrix}
0 & 1\\
0 & 0
\end{pmatrix}\Phi^{-1}
\end{align*}
and
\begin{align*}
\varphi_v&=\frac{1}{2}(-\varphi_x+\varphi_y)\\
         &=\frac{1}{2}e^{\omega/2}{\rm Ad}(-\i+\j^{'})\\
         &=e^{\omega/2}\Phi\begin{pmatrix}
0 & 0\\
1 & 0
\end{pmatrix}\Phi^{-1}.
\end{align*}
It follows from Gau{\ss}-Weingarten equations \eqref{eq:g-w} that the
coordinate frame $\Phi$ satisfies the \emph{Lax equations}:
\begin{equation}
\label{eq:lax}
\Phi_u=\Phi U,\ \Phi_v=\Phi V,
\end{equation}
where
$$U=\begin{pmatrix}
\frac{\omega_u}{4} & \frac{H}{2}e^{\omega/2}\\
-Qe^{-\omega/2} & -\frac{\omega_u}{4}
\end{pmatrix}\ {\rm and}\ V=\begin{pmatrix}
-\frac{\omega_v}{4} & Re^{-\omega/2}\\
-\frac{H}{2}e^{\omega/2} & \frac{\omega_v}{4}
\end{pmatrix}.$$
\section{Weierstra{\ss} Representation Formula of Timelike Surfaces of
  Constant Mean Curvature in Minkowski $3$-Space}

In this section, we derive the Weierstra{\ss} representation formula of
timelike surfaces of constant mean curvature in Minkowski $3$-space $\E^3_1$.

Let $M$ be a connected orientable $2$-manifold with globally defined null
coordinates system $(u,v)$. Let $\Phi: M\longrightarrow{\rm
SL}(2,\R)$ be a solution to the Lax equations \eqref{eq:lax}.

Let $\hat{\Phi}:=e^{\frac{\omega}{4}}\Phi$. Then
$\det\hat{\Phi}=e^{\frac{\omega}{2}}$.
Now, we have the following Lax equations in terms of $\hat{\Phi}$:
\begin{equation}
\label{eq:lax2}
\frac{\d\hat{\Phi}}{\d u}=\hat{\Phi}\begin{pmatrix}
\frac{\omega_u}{2} & \frac{H}{2}e^{\omega/2}\\
-Qe^{-\omega/2} & 0
\end{pmatrix},\
\frac{\d\hat{\Phi}}{\d v}=\hat{\Phi}\begin{pmatrix}
0 & Re^{-\omega/2}\\
-\frac{H}{2}e^{\omega/2} & \frac{\omega_v}{2}
\end{pmatrix}.
\end{equation}
It follows from the Lax equations \eqref{eq:lax2} that
\begin{proposition}
The conformal frame $\hat\Phi$ satisfies the \emph{Dirac equation}
\begin{equation}
\label{eq:dirac} e^{-\omega/2}\begin{pmatrix} 0 & \partial_u\\
-\partial_v & 0
\end{pmatrix}\hat\Phi^T=\frac{1}{2}H\hat\Phi^T.
\end{equation}
\end{proposition}
Conversely,
\begin{theorem}[Weierstra{\ss} Representation Formula]
Let $(s_1,-t_2)^T, (t_1,s_2)^T: M\longrightarrow\E^2_1(u,v)$ be solutions to
the Dirac equations with the potential $p\in C^\infty(M)$:
\begin{equation}
\label{eq:dirac2}
\begin{pmatrix}
 0 & \partial_u\\
-\partial_v & 0
\end{pmatrix}\begin{pmatrix}
s_1\\
-t_2
\end{pmatrix}=p\begin{pmatrix}
s_1\\
-t_2
\end{pmatrix},\ \begin{pmatrix}
 0 & \partial_u\\
-\partial_v & 0
\end{pmatrix}\begin{pmatrix}
t_1\\
s_2
\end{pmatrix}=p\begin{pmatrix}
t_1\\
s_2
\end{pmatrix}.
\end{equation}
Then $\hat\Phi:=\begin{pmatrix}
s_1 & -t_2\\
t_1 & s_2
\end{pmatrix}: M\longrightarrow\H^{'}_*$ is a conformal frame of the
conformally immersed timelike surface
$\varphi=(\varphi_1,\varphi_2,\varphi_2): M\longrightarrow\E^3_1$, where
\begin{equation}
\label{eq:weierstrass}
\begin{aligned}
\varphi_1&=\frac{1}{2}\int(s_1^2+t_1^2)du-(s_2^2+t_2^2)dv,\\
\varphi_2&=\frac{1}{2}\int(s_1^2-t_1^2)du+(s_2^2-t_2^2)dv,\\
\varphi_3&=\int (-s_1t_1du-s_2t_2dv).
\end{aligned}
\end{equation}
The metric $ds^2_\varphi$ and the mean curvature $H$ of $\varphi$ are given by
$$ds^2_\varphi=(s_1s_2+t_1t_2)^2dudv,\ H=2pe^{-\omega/2}.$$
\end{theorem}
\begin{proof}
The partial derivatives $\varphi_u$ and $\varphi_v$ are given in
terms of $\hat\Phi$ by
$$\varphi_u=e^{\frac{\omega}{2}}\hat{\Phi}\begin{pmatrix}
0 & 1\\
0 & 0
\end{pmatrix}{\hat\Phi}^{-1},\
\varphi_v=e^{\frac{\omega}{2}}\hat{\Phi}\begin{pmatrix}
0 & 0\\
1 & 0
\end{pmatrix}{\hat\Phi}^{-1}.$$

Let $\hat{\Phi}=\begin{pmatrix}
s_1 & -t_2\\
t_1 & s_2
\end{pmatrix}$ with $\det\hat{\Phi}=e^{\omega/2}$. Then
\begin{align*}d\varphi&=\varphi_udu+\varphi_vdv\\
        &=e^{\frac{\omega}{2}}\hat{\Phi}\left\{\begin{pmatrix}
0 & 1\\
0 & 0
\end{pmatrix}du+\begin{pmatrix}
0 & 0\\
1 & 0
\end{pmatrix}dv\right\}\hat{\Phi}^{-1}\\
        &=e^{\frac{\omega}{2}}\hat{\Phi}\begin{pmatrix}
0 & du\\
dv & 0
\end{pmatrix}\hat{\Phi}^{-1}\\
        &=\begin{pmatrix}
-s_1t_1du-s_2t_2dv & s_1^2du-t_2^2dv\\
-t_1^2du+s_2^2dv & s_1t_1du+s_2t_2dv
\end{pmatrix}.
\end{align*}
Now, $\varphi$ can be retrieved by integrating $d\varphi$:
\begin{align*}
\varphi&=\int d\varphi\\
       &=\int\begin{pmatrix}
-s_1t_1du-s_2t_2dv & s_1^2du-t_2^2dv\\
-t_1^2du+s_2^2dv & s_1t_1du+s_2t_2dv
\end{pmatrix}.
\end{align*}
On the other hand, by the identification \eqref{eq:matrix},
$$\varphi=(\varphi_1,\varphi_2,\varphi_3)\cong\begin{pmatrix}
\varphi_3 & \varphi_1+\varphi_2\\
-\varphi_1+\varphi_2 & -\varphi_3
\end{pmatrix}.$$
Hence, we obtain
\begin{align*}
\varphi_1&=\frac{1}{2}\int(s_1^2+t_1^2)du-(s_2^2+t_2^2)dv,\\
\varphi_2&=\frac{1}{2}\int(s_1^2-t_1^2)du+(s_2^2-t_2^2)dv,\\
\varphi_3&=\int (-s_1t_1du-s_2t_2dv).
\end{align*}
\end{proof}
\section{The Gau{\ss} Map and Weierstra{\ss} Representation Formula}
In the previous section, we obtained Weierstra{\ss} representation formula for
timelike surfaces in Minkowski $3$-space. In this section, we study the
relationship between the data $s_1,s_2,s_1,t_2$ and the Gau{\ss} map of
timelike surface which is given by the Weierstra{\ss} representation formula.

The pseudosphere of radius $1$ in Minkowski $3$-space is the hyperquadric
$$S^2_1=\{(x_1,x_2,x_3)\in\E^3_1: -x_1^2+x_2^2+x_3^2=1\}$$
of constant Gau{\ss}ian curvature $1$.

Let $\varphi: M\longrightarrow\E^3_1$ be an orientable timelike surface and
$N$ a unit normal vector field to $\varphi$. This unit normal vector field
$p\in M\longmapsto N(p)\in S^2_1$
is defined to be the \emph{Gau{\ss} map} of $\varphi$.

The Ad-action of ${\rm SL}(2,\R)$ on $S^2_1$ is transitive and isometric. The
sotropy subgroup at $\k^{'}$ is ${\rm SO}(1,1)$. Thus, $S^2_1$ is identified
with the homogeneous space
$${\rm SL}(2,\R)/{\rm
    SO}(1,1)=\left\{h\k^{'}h^{-1}: h\in{\rm SL}(2,\R)\right\}.$$
Let $\wp_{\mathcal N}: S^2_1\setminus\{x_3=1\}\longrightarrow\E^2_1\setminus
  H_0^1$ be the stereographic projection from the north pole
  ${\mathcal N}=(0,0,1)$. Here, $H_0^1=\{(x_1,x_2)\in\E^2_1:
  -x_1^2+x_2^2=-1\}$. Then
\begin{align*}
\wp_{\mathcal N}(x_1,x_2,x_3)&=\left(\frac{x_1}{1-x_3},\frac{x_2}{1-x_3}\right)\\
                  &\cong\left(\frac{x_1+x_2}{1-x_3},\frac{-x_1+x_2}{1-x_3}\right)\
                  \mbox{\rm in null coordinate system}\ (u,v).
\end{align*}
The inverse stereographic projection $\wp_{\mathcal N}^{-1}: \E^2_1\setminus
H_0^1\longrightarrow S^2_1\setminus\{x_3=1\}$ is given by
$$\wp_{\mathcal
  N}^{-1}(x,y)=\left(\frac{2x}{1-x^2+y^2},\frac{2y}{1-x^2+y^2},\frac{-1-x^2+y^2}{1-x^2+y^2}\right).$$
In terms of null coordinates $(u,v)$,
\begin{equation}
\label{eq:invproj}
\wp_{\mathcal N}^{-1}(x,y)=\wp_{\mathcal N}^{-1}\left(\frac{u-v}{2},\frac{u+v}{2}\right)=\left(\frac{u-v}{1+uv},\frac{u+v}{1+uv},\frac{-1+uv}{1+uv}\right).
\end{equation}

The Gau{\ss} map can be written as $N=
h\begin{pmatrix}
1 & 0\\
0 & -1
\end{pmatrix}h^{-1}\in S^2_1$ for some $h=\begin{pmatrix}
p_1 & -q_2\\
q_1 & p_2
\end{pmatrix}\in{\rm SL}(2,\R)$.
That is, $N=\begin{pmatrix}
p_1p_2-q_1q_2 & 2p_1q_2\\
2p_2q_1 & -p_1p_2+q_1q_2
\end{pmatrix}$ and its projected image via $\wp_{\mathcal N}$ is
$\wp_{\mathcal N}(N)\cong\left(\frac{p_1}{q_1},\frac{p_2}{q_2}\right)$ in null
coordinate system $(u,v)$.

Now, $\hat{\Phi}=\begin{pmatrix}
s_1 & -t_2\\
t_1 & s_2
\end{pmatrix}=e^{\frac{\omega}{2}}\begin{pmatrix}
p_1 & -q_2\\
q_1 & p_2
\end{pmatrix}$ and $\begin{pmatrix}
p_1 & -q_2\\
q_1 & p_2
\end{pmatrix}\in{\rm SL}(2,\R)$. Then,
$$q:=\frac{p_1}{q_1}=\frac{s_1}{t_1}\ {\rm and}\
r:=\frac{p_2}{q_2}=\frac{s_2}{t_2}.$$
Therefore, the Weierstra{\ss} representation formula \eqref{eq:weierstrass}
can be written in terms of the projected Gau{\ss} map $(q,r)$ as
\begin{equation}
\label{eq:weierstrass2}
\begin{aligned}
\varphi_1&=\frac{1}{2}\int
(1+q^2)f(u)du-(1+r^2)g(v)dv,\\
\varphi_2&=-\frac{1}{2}\int
(1-q^2)f(u)du+(1-r^2)g(v)dv,\\
\varphi_3&=-\int qf(u)du+rg(v)dv,
\end{aligned}
\end{equation}
where $f(u)=t_1^2$ and $g(v)=t_2^2$. The metric is given by
\begin{equation}
\label{eq:metric}
ds^2_\varphi=(1+qr)^2f(u)g(v)dudv.
\end{equation}

Using the Lax equations \eqref{eq:lax2}, we compute
\begin{proposition}
The projected Gau{\ss} map
$\wp_{\mathcal N}\circ N$ satisfy the following equations:
\begin{align}
\label{eq:gmap1}
q_u&=\frac{Q}{f(u)},\\
\label{eq:gmap2}
q_v&=\frac{He^{\omega}}{2f(u)}=\frac{H}{2}(1+qr)^2g(v),\\
\label{eq:gmap3}
r_u&=\frac{He^{\omega}}{2g(v)}=\frac{H}{2}(1+qr)^2f(u),\\
\label{eq:gmap4}
r_v&=\frac{R}{g(v)}.
\end{align}
\end{proposition}
\begin{remark}
In \cite{I-T}, the authors studied the \emph{normalized
Weierstra{\ss} formula} for timelike minimal surfaces via \emph{loop group
method}. In this case, 
$f(u)=g(v)=1$ and so, by equations \eqref{eq:gmap1} and
\eqref{eq:gmap4}, $q$ and $r$ are interpreted as the primitive
functions of the coefficients $Q$ and $R$, resp. of the Hopf
differential \eqref{eq:hopf}. Originally, the Weierstra{\ss} formula
\eqref{eq:weierstrass2} was 
obtained by M. A. Magid in \cite{magid}, 
however the geometric meaning of the data $(q,r)$ is not clarified.
In \cite{I-T}, the data $(q,r)$ are retrieved from the
\emph{normalized potential} in their construction.
\end{remark}
The corollaries \ref{cor:min} and \ref{cor:umb} immediately follow
from the equations \eqref{eq:gmap1}-\eqref{eq:gmap4}.
\begin{corollary}
\label{cor:min}
A Lorentz surface $\varphi: M\longrightarrow\E^2_1$
is minimal (i.e., $H=0$) if and only if $q$ is Lorentz holomorphic
and $r$ is Lorentz antiholomorphic.
\end{corollary}
\begin{corollary}
\label{cor:umb}
A Lorentz surface $\varphi: M\longrightarrow\E^2_1$
is totally umbilic if and only if $q$ is Lorentz antiholomorphic and
$r$ is Lorentz holomorphic.
\end{corollary}
\begin{corollary}
\label{cor:tum}
A totally umbilic minimal Lorentz surface in $\E^3_1$ is part of timelike
plane in $\E^3_1$.
\end{corollary}
\begin{proof}
By the equation \eqref{eq:sff}, a totally umbilic minimal Lorentz
surface has $I\!I=0$, i.e., it is totally geodesic.
\end{proof}
\begin{theorem}[Weierstra{\ss} Representation Formula for Timelike Minimal
    Surfaces in $\E^3_1$]
Let $q,\ r: M\longrightarrow\E^2_1$ be Lorentz holomorphic and
Lorentz antiholomorphic maps, resp. Then
\begin{equation}
\label{eq:weierstrass3}
\begin{aligned}
\varphi_u&=\left(\frac{1}{2}(1+q^2),-\frac{1}{2}(1-q^2),-q\right)f(u),\\
\varphi_v&=\left(-\frac{1}{2}(1+r^2),-\frac{1}{2}(1-r^2),-r\right)g(v)
\end{aligned}
\end{equation}
define a timelike minimal surface $\varphi: M\longrightarrow\E^2_1$.
Here, $f(u)$ and $g(v)$ are Lorentz holomorphic and Lorentz
antiholomorphic maps. The metric of $\varphi$ is given by
$ds_\varphi^2=(1+qr)^2f(u)g(v)dudv$. Conversley, any timelike minimal surface
can be represented by \eqref{eq:weierstrass3} up to translations.
\end{theorem}
\begin{proof}
The first statement follows immediately from \eqref{eq:weierstrass2},
\eqref{eq:metric}, \eqref{eq:gmap2}, and \eqref{eq:gmap3}.

Let $\varphi: M\longrightarrow\E^3_1$ be a timelike minimal
surface parametrized by null coordinates $(u,v)$. Define two vector valued
functions $\xi=(\xi_0,\xi_1,\xi_2)$ and $\eta=(\eta_0,\eta_1,\eta_2)$ by 
$$\xi(u):=\varphi_u,\ \eta(v):=\varphi_v,$$
that is,
$$X(u)=\int_0^u\xi(u)du,\ Y(v)=\int_0^v\eta(v)dv$$
for the timelike minimal surface $\varphi(u,v)=X(u)+Y(v)$. Since $\varphi$ is
a solution to the wave equation $\varphi_{uv}=0$,
$\xi$ and $\eta$ are a Lorentz holomorphic and a Lorentz anti-holomorphic null
curves, resp., in $\E^3_1$. Hence, they satisfy
$$-\xi_0^2+\xi_1^2+\xi_2^2=-\eta_0^2+\eta_1^2+\eta_2^2=0.$$
Define the functions $q(u), f(u), r(v)$ and $g(v)$ by
\begin{align*}
-\xi_0+\xi_1&=-f,\ \xi_2=-qf,\\
\eta_0+\eta_1&=-g,\ \eta_2=-rg.
\end{align*}
Then we obtain
\begin{align*}
\xi(u)&=\left(\frac{1}{2}(1+q(u)^2),-\frac{1}{2}(1-q(u)^2),-q(u)\right)f(u),\\
\eta(v)&=\left(-\frac{1}{2}(1+r(v)^2),-\frac{1}{2}(1-r(v)^2),-r(v)\right)g(v). 
\end{align*}
Thus, the given
timelike minimal surface is represented by 
\begin{equation}
\begin{aligned}
\varphi(u,v)&=\int_0^u\left(\frac{1}{2}(1+q(u)^2),-\frac{1}{2}(1-q(u)^2),-q(u)\right)f(u)du\\
       &+\int_0^v\left(-\frac{1}{2}(1+r(v)^2),-\frac{1}{2}(1-r(v)^2),-r(v)\right)g(v)dv
\end{aligned}
\end{equation}
up to translations.
\end{proof}
\begin{remark}
The maps $q$ and $r$ coincide with the first and the second
component maps of the projected Gau{\ss} map $\wp_{\mathcal N}\circ
N$ of the timelike minimal surface $\varphi$.
\end{remark}
\begin{example}[Lorentz Enneper surfaces]

Let $q(u)=\varepsilon u=\pm u$, $f(u)=1$, $r(v)=v$, $g(v)=1$. Then we obtain
the following timelike minimal immersion  
$$\varphi^{(\varepsilon)}(u,v)=X(u)+Y(v),$$
where
\begin{align*}
X(u)&=\frac{1}{2}\left(u+\frac{u^3}{3},-u+\frac{u^3}{3},\mp u^2\right),\\
Y(v)&=\frac{1}{2}\left(-v-\frac{v^3}{3},-v+\frac{v^3}{3},-v^2\right).
\end{align*}
The timelike minimal surface $\varphi^{(\varepsilon)}$ is called \emph{Lorentz
  Enneper surface}. 

The metric of $\varphi^{(\varepsilon)}$ is given by 
$$I=(1+\varepsilon uv)^2dudv.$$
The Hopf differential and the Gau{\ss}ian curvature of $\varphi^{(\varepsilon)}$
are 
$$\varepsilon du^2+dv^2,\ K=-4\varepsilon(1+\varepsilon uv)^{-4}.$$
The surface $\varphi^{(-1)}$ has two imaginary principal curvatures, while
$\varphi^{(1)}$ has real distinct principal curvatures.

\begin{figure}[ht]
\centering
\mbox{\subfigure[]{\epsfig{figure=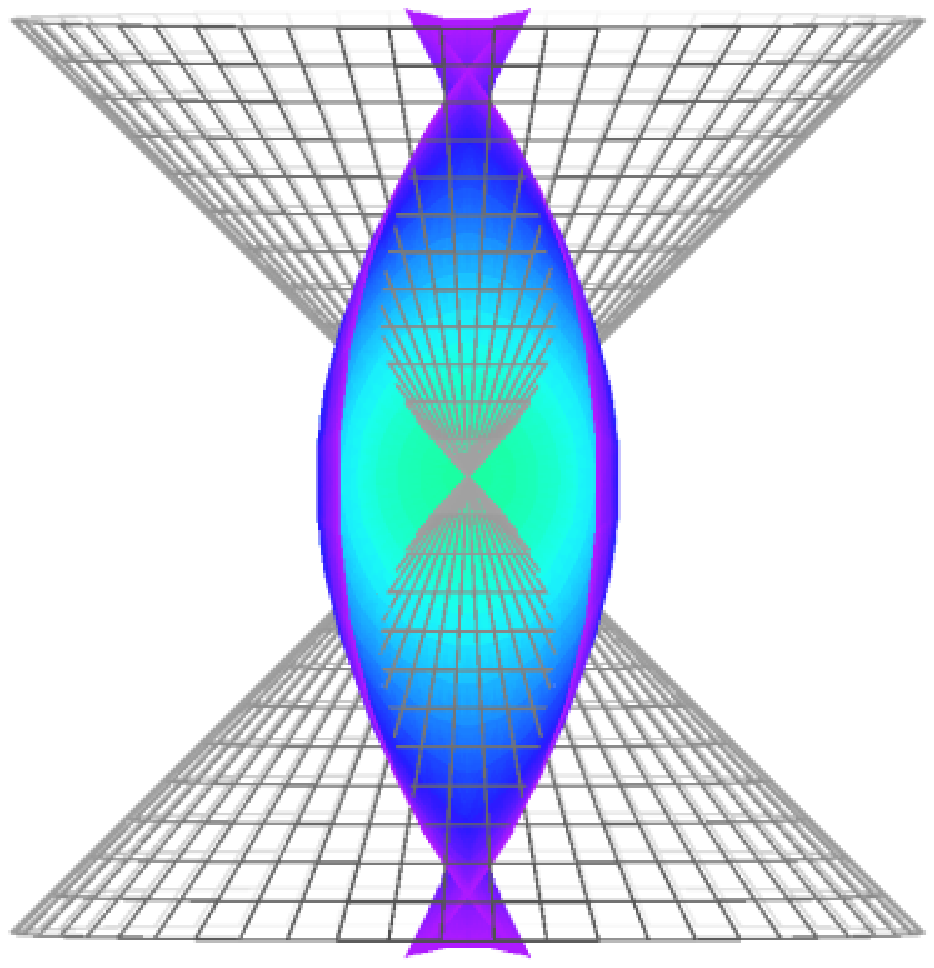,width=.40\textwidth}}\quad
    \subfigure[]{\epsfig{figure=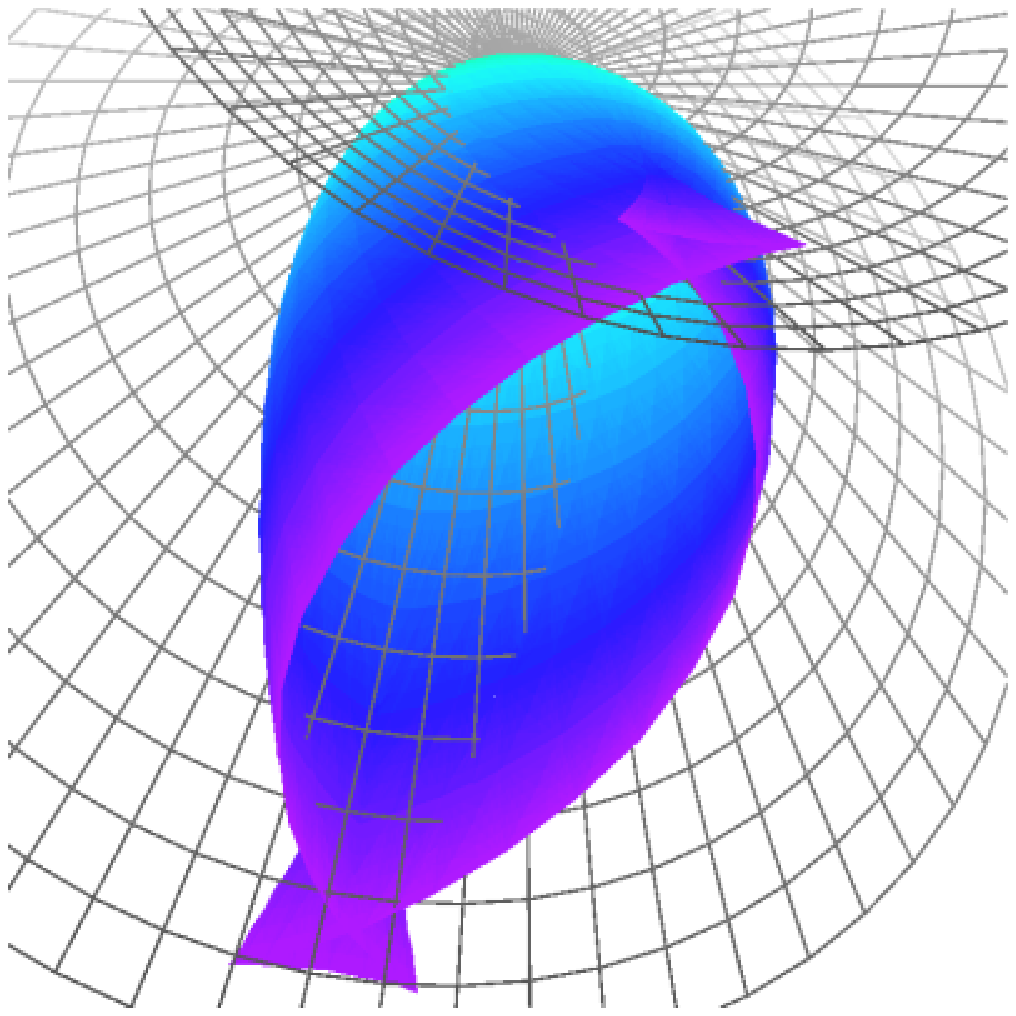,width=.40\textwidth}}}
\mbox{\subfigure[]{\epsfig{figure=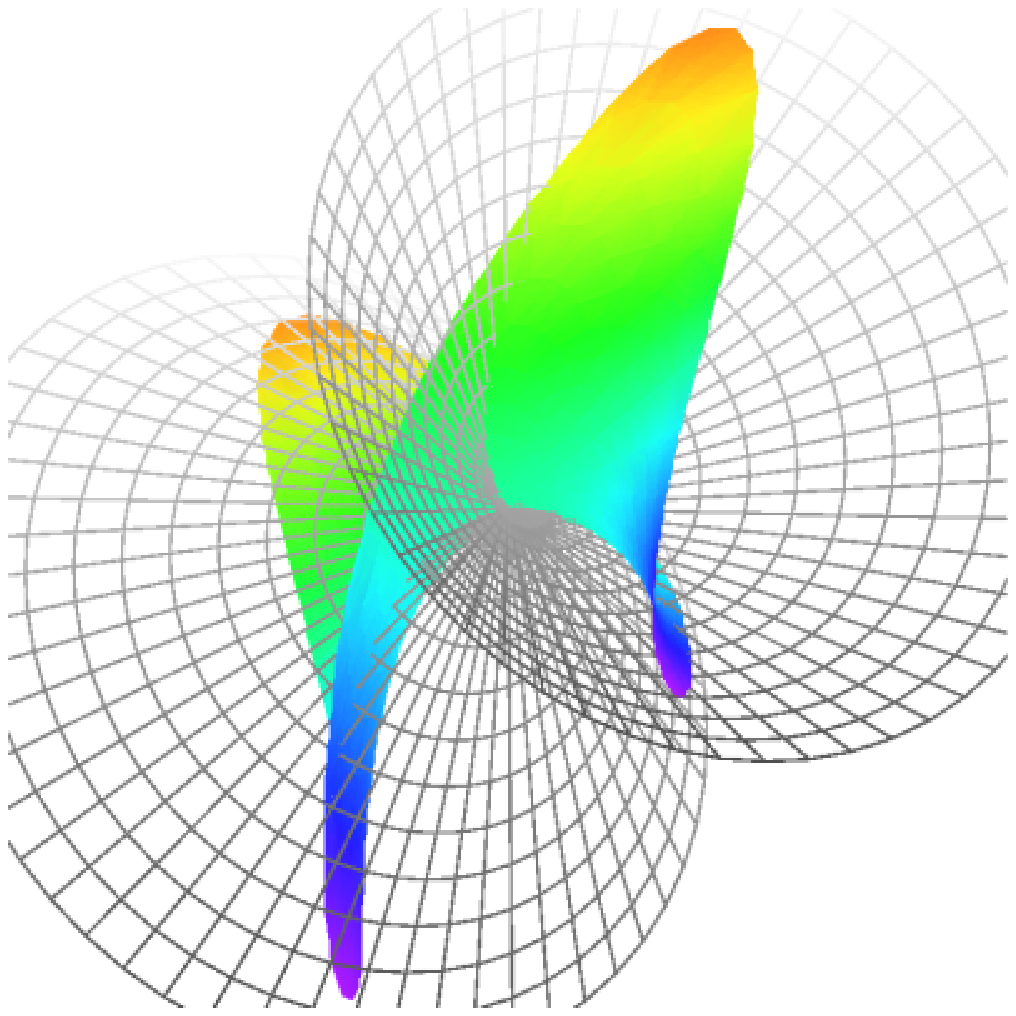,width=.40\textwidth}}\quad
    \subfigure[]{\epsfig{figure=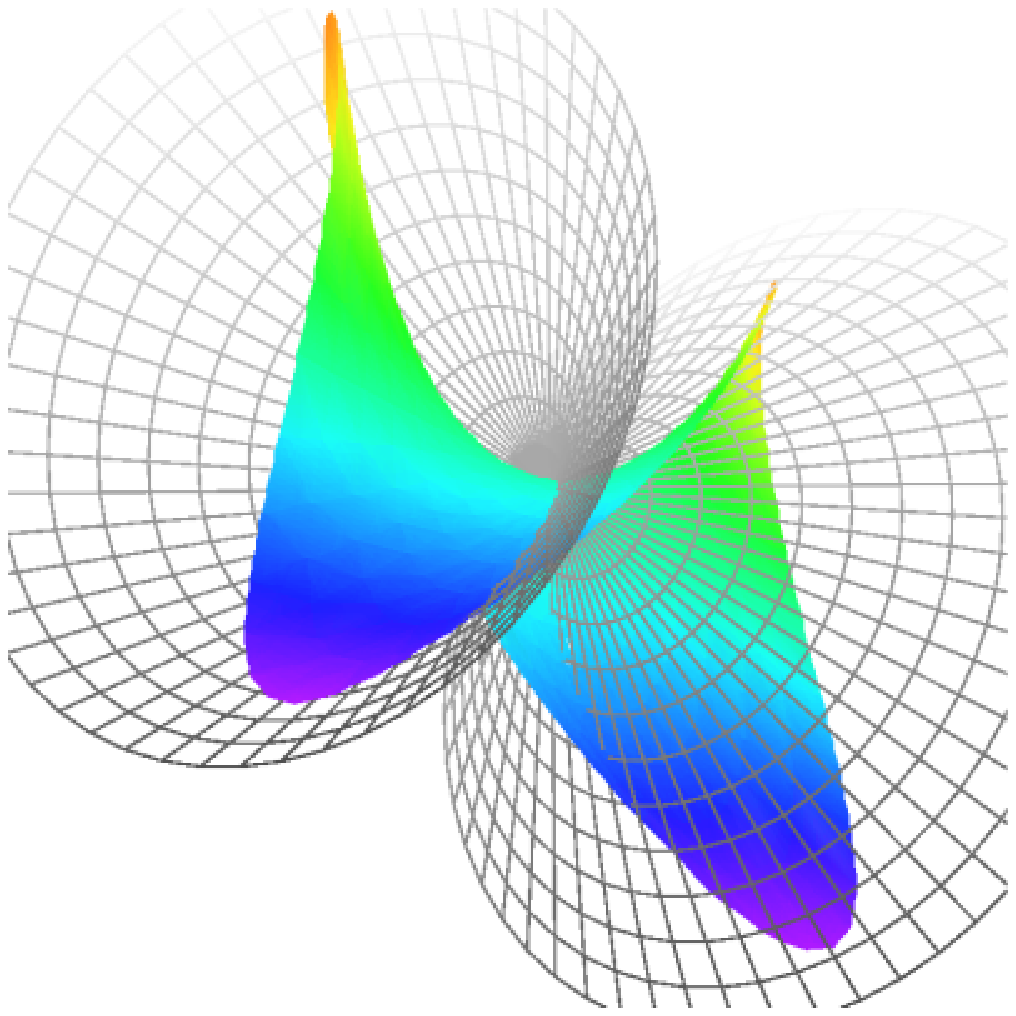,width=.40\textwidth}}}
\caption{Lorentz Enneper surfaces $\varphi^{(1)}$ (a), (b) and
$\varphi^{(-1)}$ (c), (d)}
\end{figure}
\end{example}
\begin{definition}
Let $\varphi(u,v)=X(u)+Y(v)$ be a timelike minimal surface. Then clearly,
$\hat\varphi(u,v):=X(u)-Y(v)$ is also a timelike minimal surface. The timelike
minimal surface $\hat\varphi$ is called the \emph{conjugate timelike minimal
  surface} of $\varphi$.
\end{definition} 
\begin{example}[Lorentz catenoid and Lorentz helicoid with a spacelike axis]

\emph{Lorentz catenoids} are timelike minimal surfaces of revolution with a
spacelike 
axis or timelike axis. Lorentz catenoid with a spacelike axis 
$\varphi(u,v)=X(u)+Y(v)$ (Figure \ref{fig:tlmincatenoid1}) can be obtained
by the Weierstra{\ss} formula \eqref{eq:weierstrass3} with data
$q(u)=-e^u$, $f(u)=-e^{-u}$, $r(v)=e^{-v}$, $g(v)=-e^v$, where
$$
X(u)=(-\sinh u,-\cosh u,-u),\ Y(v)=(\sinh v,\cosh v,v).$$

The conjugate surface $\hat\varphi=X(u)-Y(v)$ is called \emph{Lorentz
  helicoid} with a spacelike axis (Figure \ref{fig:tlminhelicoid1}).

\begin{figure}[ht]
\centering
\mbox{\subfigure[]{\epsfig{figure=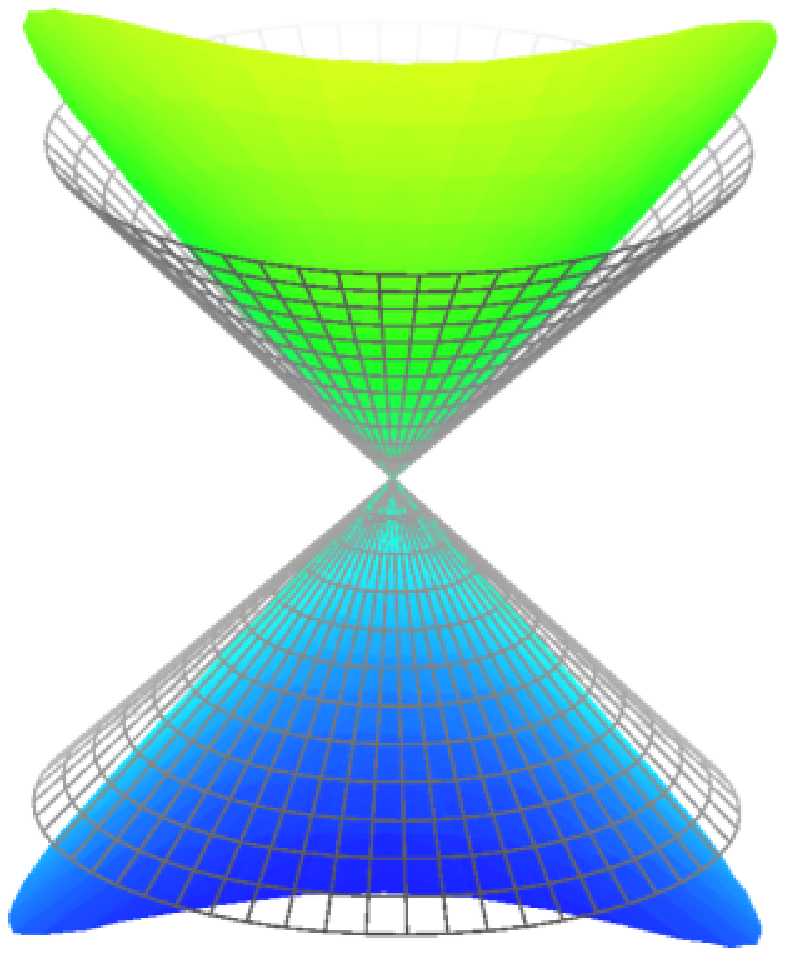,width=.40\textwidth}\label{fig:tlmincatenoid1}}\quad
    \subfigure[]{\epsfig{figure=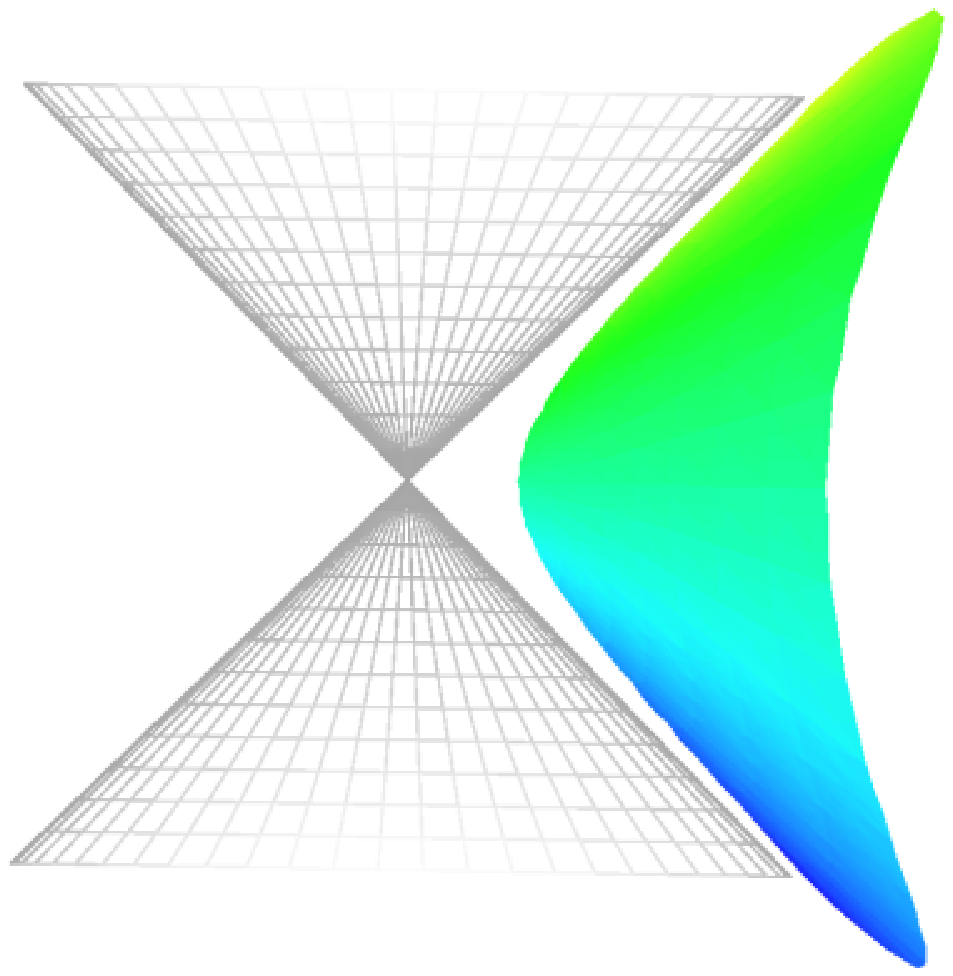,width=.40\textwidth}\label{fig:tlminhelicoid1}}}   
\caption{Lorentz catenoid and Lorentz helicoid with a spacelike axis}
\end{figure}  

The ordered pair $(q,r)=(-e^u,e^{-v})$ is the Gau{\ss} map projected
in $\E^2_1(u,v)$. Note that Lorentz catenoid and Lorentz helicoid with a
spacelike axis share the same Gau{\ss} map analogously to the Euclidean case.
\end{example}
\begin{example}[Lorentz catenoid and Lorentz helicoid with timelike axis]

Lorentz catenoid with timelike axis 
$\varphi(u,v)=X(u)+Y(v)$ (Figure \ref{fig:tlmincatenoid2}) can be obtained
by the Weierstra{\ss} formula \eqref{eq:weierstrass3} with data 
$q(u)=\frac{\sin u}{-1+\cos u}$, $f(u)=-1+\cos u$, $r(v)=\frac{\sin
  v}{1+\cos v}$, $g(v)=-(1+\cos v)$, where
$$
X(u)=(-u,-\sin u,\cos u),\ Y(v)=(v,\sin v,-\cos v).$$

The conjugate surface $\hat\varphi=X(u)-Y(v)$ is called \emph{Lorentz
 helicoid} with timelike axis (Figure \ref{fig:tlminhelicoid2}).

\begin{figure}[ht]
\centering
\mbox{\subfigure[]{\epsfig{figure=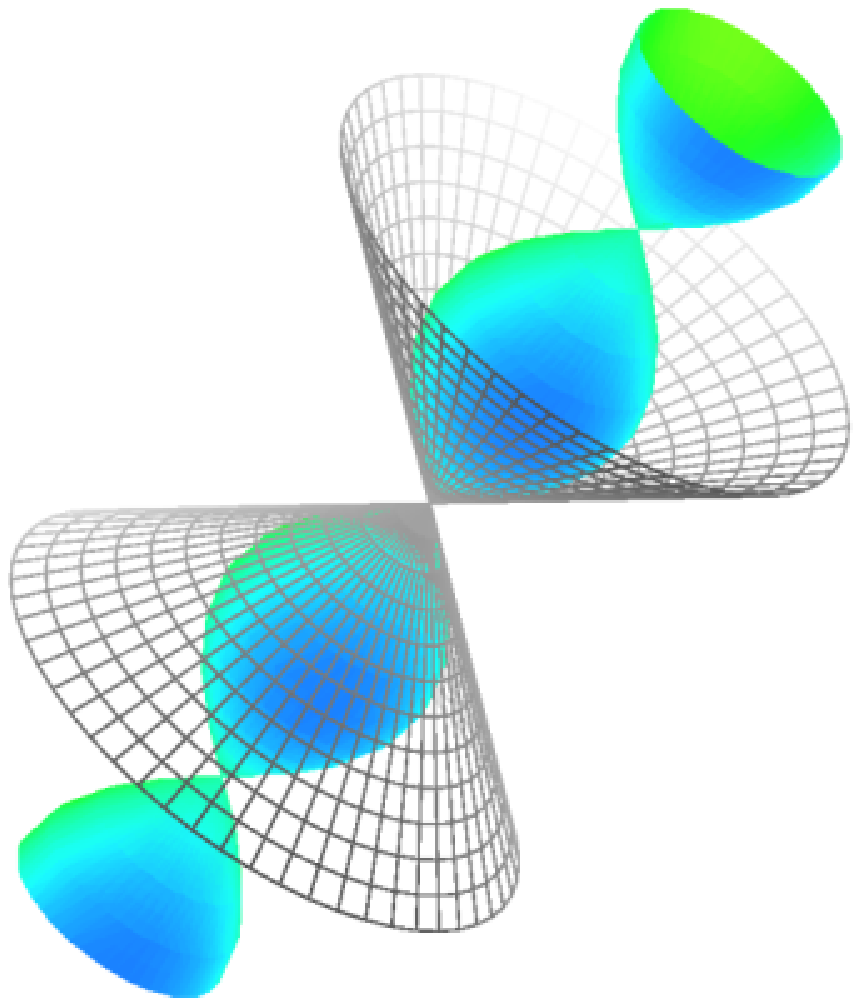,width=.40\textwidth}\label{fig:tlmincatenoid2}}\quad
    \subfigure[]{\epsfig{figure=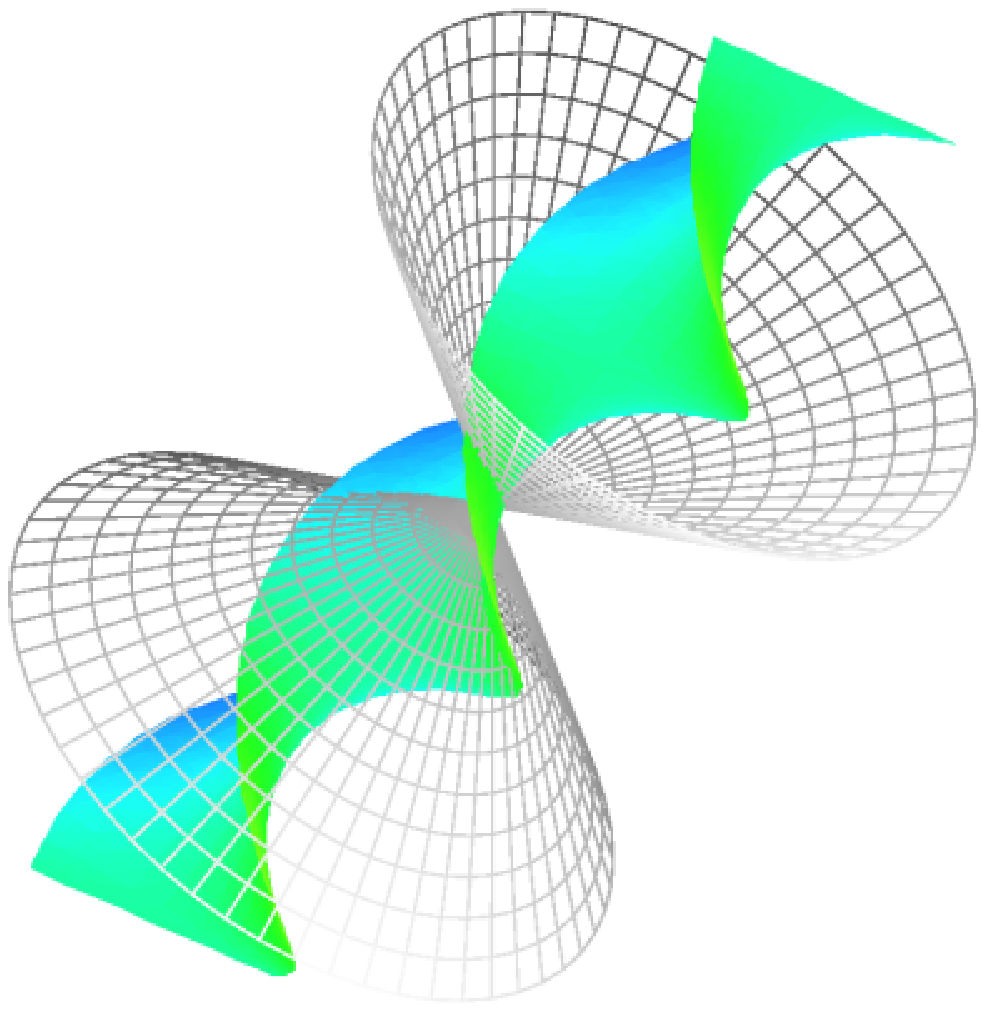,width=.40\textwidth}\label{fig:tlminhelicoid2}}}   
\caption{Lorentz catenoid and Lorentz helicoid with timelike axis}
\end{figure}

The ordered pair $(q,r)=\left(\frac{\sin u}{-1+\cos u},\frac{\sin
  v}{1+\cos v}\right)$ is the Gau{\ss} map projected
in $\E^2_1(u,v)$. Note that timelike catenoid and timelike helicoid with
  timelike axis share the same Gau{\ss} map analogously to the Euclidean
  case. 
\end{example}

It is known that a totally umbilic Lorentz surface in $\E^3_1$ is
congruent to an open part of the pseudosphere $S^2_1$ or timelike
plane in $\E^3_1$. This property can be proved directly from the
Weierstra{\ss} formula \eqref{eq:weierstrass2}.
\begin{corollary}
A totally umbilic timelike surface in $\E^3_1$ is part of timelike
plane or a pseudosphere in $\E^3_1$. In particular, a totally
umbilic timelike surface with positive Gau{\ss}ian curvature $K$ is
part of a pseudosphere in $\E^3_1$ of radius $\frac{1}{\sqrt{K}}$.
\end{corollary}
\begin{proof}
Let $M$ be a simply connected domain that contains the origin
$(0,0)$. Let $\varphi=(\varphi_1,\varphi_2,\varphi_3):
M\longrightarrow\E^3_1$ be a totally umbilic timelike surface in
$\E^3_1$. If $\varphi$ is minimal,by Corollary \ref{cor:tum} it is
part of timelike plane. So, we assume that $\varphi$ is not minimal,
i.e., $H\ne 0$.

By the Weierstr{\ss} formula \eqref{eq:weierstrass2},
\begin{align*}
\varphi_1&=\int\frac{1}{2}(1+q^2)f(u)du-\frac{1}{2}(1+r^2)g(v)dv\\
         &=\int\frac{1+q^2}{H(1+qr)^2}dr-\frac{1+r^2}{H(1+qr)^2}dq.
\end{align*}
Since $\dis\frac{\d\varphi_1}{\d r}=\frac{1+q^2}{H(1+qr)^2}$,
\begin{align*}
\varphi_1&=\int_0^r\frac{1+q^2}{H(1+qr)^2}dr+\psi(q)\\
         &=-\frac{1+q^2}{H(1+qr)q}+\psi(q).
\end{align*}
So,
$$\frac{\d\varphi_1}{\d
  q}=\frac{-q^2+2qr+1}{(1+qr)^2}+\psi^{'}(q).$$
This must be the same as $-\dis\frac{1+r^2}{H(1+qr)^2}$; hence
  $\psi^{'}(q)=-\dis\frac{1}{Hq^2}$, that is,
  $\psi(q)=\dis\frac{q}{Hq}+c_1,$ where $c_1$ is a constant.
Thus, $\varphi_1$ is given by
$$\varphi_1=-\frac{q-r}{H(1+qr)}+c_1.$$
Similarly, we compute
\begin{align*}
\varphi_2&=-\frac{q+r}{H(1+qr)}+c_2,\\
\varphi_3&=\frac{1-qr}{H(1+qr)}+c_3,
\end{align*}
where $c_2$ and $c_3$ are constants.
The parametric equation
$$\varphi=\left(-\frac{q-r}{H(1+qr)}+c_1,-\frac{q+r}{H(1+qr)}+c_2,\frac{1-qr}{H(1+qr)}+c_3\right)$$
shows that $\varphi$ is part of the pseudosphere in $\E^3_1$
centered at $(c_1,c_2,c_3)$ with radius $\frac{1}{H}$. (See the
equation \eqref{eq:invproj}.) Since $\varphi$ is totally umbilic,
$Q=R=0$. So, by the equation \eqref{eq:gauss2}, $H^2=K$. If the
Gau{\ss}ian curvature $K$ is positive, then $H=\sqrt{K}$. This
completes the proof.
\end{proof}
\section{Appendix: Timelike minimal surfaces and bosonic Nambu-Goto strings in
  Minkowski spacetime} 
By the equation \eqref{eq:meancurv2},  a Lorentz surface $\varphi:
M\longrightarrow\E^3_1$ is minimal if and only if it is a solution
to the homogeneous \emph{wave equation} $\varphi_{uv}=0$ (or
equivalently $\Box\varphi=0$). Hence, as seen in
\eqref{eq:weierstrass2}, timelike minimal surfaces in Minkowski
$(2+1)$-spacetime can be retrieved by integrating a pair of Lorentz
holomorphic and Lorentz antiholomorphic null curves in spacetime,
i.e., curves in the light cone ${\mathcal
N}=\{(x_1,x_2,x_3)\in\E^3_1: -x_1^2+x_2^2+x_3^2=0\}$. Note that
Lorentz holomorphic and Lorentz antiholomorphic null curves in
spacetime are the trajectories of massless particles (\emph{bosonic}
particles) in spacetime. Therefore, there appears to be some
relationship  between timelike minimal surfaces and bosonic
particles in spacetime. In string theory, elementary particles are
considered to be tiny vibrating strings in spacetime. A string
evolves in time while sweeping a surface, the so-called
\emph{worldsheet}, in spacetime. Hence, string worldsheets are
timelike surfaces. Moreover, it can be shown that string worldsheets
are indeed timelike minimal surfaces in spacetime. The following
discussion is not restricted in $(2+1)$-dimensions. So, one can
assume the standard $(3+1)$-spacetime or any higher
(D+1)-dimensional spacetime. (Here, we consider only classical strings. In
case of \emph{quantized} strings, our universe is, without
supersymmetry, a $26$-dimensional spacetime and it is a $10$-dimensional
spacetime with supersymmetry.)

In order to be consistent with physicists' common notations, we
consider string worldsheets parametrized by $\tau$ and $\sigma$,
where $\tau$ is time parameter in $(1+1)$-spacetime. The motion of
bosonic strings in spacetime is described by the \emph{Nambu-Goto
string action}
\begin{equation}
\label{eq:action} {\mathcal S}=-T\int(-\det h_{ab})^{1/2}d\tau
d\sigma,
\end{equation}
where $T$ is tension and $h_{ab}$ is the metric tensor of the
worldsheet $\varphi: M\longrightarrow\E^D_1$. Using the
\emph{Einstein's convention}, $h_{ab}$ is given by
$$h_{ab}=\partial_a\varphi^\mu\partial_b\varphi^\nu\eta_{\mu\nu},$$
where $\eta_{\mu\nu}$ is the metric tensor of the flat Minkowski
$(D+1)$-spacetime with signature $(-,+,\cdots,+)$. Clearly,
$dA:=(-\det h_{ab})^{1/2}d\tau d\sigma$ is the area element of the
string worldsheet $\varphi$.

Let us denote $\dot\varphi:=\frac{\partial\varphi}{\partial\tau}$
and $\varphi^{'}:=\frac{\partial\varphi}{\partial\sigma}$. The
\emph{Lagrangian} ${\mathcal L}$ of the string motion is
\begin{equation}
\label{eq:lagrange} {\mathcal
L}(\dot\varphi,\varphi^{'};\sigma,\tau)=-T(-\det h_{ab})^{1/2}.
\end{equation}
By variational principle, $\delta{\mathcal S}=0$, subject to the
condition that the initial and final configurations of the string
are kept fixed, implies the \emph{Euler-Lagrange} equation for the
string action \eqref{eq:action} is
\begin{equation}
\label{eq:E-L} \frac{\partial}{\partial\tau}\frac{\partial{\mathcal
L}}{\partial\dot\varphi^\mu}+\frac{\partial}{\partial\sigma}\frac{\partial{\mathcal
L}}{\partial{\varphi^{'}}^\mu}=0
\end{equation}
with $$\frac{\partial{\mathcal
L}}{\partial\dot\varphi^\mu}=\frac{\partial{\mathcal
L}}{\partial{\varphi^{'}}^\mu}=0\ {\rm on}\ \partial M.$$
Note that
the action \eqref{eq:action} is invariant under conformal scaling of
the worldsheet metric. Physicists call it \emph{Weyl invariance} and
it is an important symmetry of the action \eqref{eq:action} along
with \emph{worldsheet reparametrizations} and
\emph{Lorentz/Poincar\'{e} symmetries}. With \emph{conformal gauge
fixing}
$$-\langle\varphi_\tau,\varphi_\tau\rangle=\langle\varphi_\sigma,\varphi_\sigma\rangle=e^\omega,\
\langle\varphi_\tau,\varphi_\sigma\rangle=0,$$ one can easily show
that the Euler-Lagrange equation \eqref{eq:E-L} is equivalent to the
homogeneous wave equation
$$\Box\varphi=-\frac{\partial^2\varphi}{\partial\tau^2}+\frac{\partial^2\varphi}{\partial\sigma^2}=0.$$
Hence, we see that string worldsheets are timelike minimal surfaces.

The Gau{\ss}ian curvature $K$ also plays
an important role in string theory, especially when results from
different orders of string perturbation theory are compared and when string
interactions are considered. In order to write a more familiar
expression in string theory, we introduce a new parameter
$\alpha^{'}$ which is defined by
$$\alpha^{'}=\frac{1}{2\pi T}.$$
This parameter $\alpha^{'}$ is known as the slope of \emph{Regge
trajectories} in physics. We refer to
\cite{kiritsis}, \cite{P-W} for physical background and details
regarding Regge trajectories. As is well-known in
variational theory, one can add more terms to action functionals as
constraints. One physically interesting extra term is the
\emph{Einstein-Hilbert action}
\begin{equation}
\label{eq:E-H} \chi:=\frac{1}{4\pi\alpha^{'}}\int_M R(-\det
h_{ab})^{1/2}d\tau d\sigma+\frac{1}{2\pi\alpha^{'}}\int_{\partial M}
kds,
\end{equation}
where $R$ is the Ricci scalar (or scalar curvature) on the
worldsheet $M$ and $k$ is the geodesic curvature on $\partial M$.
Note that the Einstein-Hilbert action \eqref{eq:E-H} is invariant
under the conformal transformation (Weyl transformation)
$h_{ab}\longrightarrow e^\omega h_{ab}$. Now, the full string action
is given by
$$
\begin{aligned}{\mathcal S}^{'}=-\frac{1}{2\pi\alpha^{'}}\int(-\det
h_{ab})^{1/2}d\tau
d\sigma&+\lambda\left\{\frac{1}{4\pi\alpha^{'}}\int_M R(-\det
h_{ab})^{1/2}d\tau d\sigma\right.\\
       &\left.+\frac{1}{2\pi\alpha^{'}}\int_{\partial M} kds\right\},
       \end{aligned}$$ where
$\lambda$ is a coupling parameter. This full string action resembles
$2$-dim gravity coupled with bosonic matter fields and the equations
of motion is given by the following Einstein's field equation:
$$R_{ab}-\frac{1}{2}h_{ab}R=T_{ab}.$$
The LHS of the Einstein's equation vanishes identically in
$2$-dimensions, so there is no dynamics associated with
\eqref{eq:E-H}. For surfaces, the Ricci scalar $R$ and the
Gau{\ss}ian curvature $K$ are related by $R=2K$, thus \eqref{eq:E-H}
can be written as
\begin{equation}
\label{eq:G-B} \chi:=\frac{1}{2\pi\alpha^{'}}\int_M
KdA+\frac{1}{2\pi\alpha^{'}}\int_{\partial M} kds.
\end{equation}
This is nothing but a constant multiple
$\left(\frac{1}{2\pi\alpha^{'}}\right)$ of the \emph{Gau{\ss}-Bonnet
term}. In string theory, physicists apply the so-called \emph{Wick
rotation} $\tau\longrightarrow i\tau$ to change Lorentzian signature
to Euclidean signature. As a result, string worldsheets turn into
Riemann surfaces with complex local coordinates $\tau\pm i\sigma$.
This procedure is required because string amplitudes are computed by
the Feynman path integral which is defined in Euclidean setting.
Once calculation is done, one retrieves the Lorentzian signature by
the opposite Wick rotation for a physical interpretation. In
Euclidean setting, by the well-known \emph{Gau{\ss}-Bonnet Theorem},
the RHS of \eqref{eq:G-B} is the same as
$\frac{1}{\alpha^{'}}\chi(M)$ where $\chi(M)$ is the Euler
characteristic of the string worldsheet $M$ as a \emph{compact Riemann
surface}. String worldsheets that are swept by closed strings are
compact orientable surfaces without boundary and so $\chi(M)=2-2g$,
where $g$ is the genus of the (Riemannian) worldsheet $M$.
Therefore, a closed string is distinguished from another by the
genus of its (Riemannian) worldsheet, which solely depends on the
topology of the worldsheet.

\end{document}